# The distribution of a linear predictor after model selection: Unconditional finite-sample distributions and asymptotic approximations

Hannes Leeb[1],*

*Yale University*

**Abstract:** We analyze the (unconditional) distribution of a linear predictor that is constructed after a data-driven model selection step in a linear regression model. First, we derive the exact finite-sample cumulative distribution function (cdf) of the linear predictor, and a simple approximation to this (complicated) cdf. We then analyze the large-sample limit behavior of these cdfs, in the fixed-parameter case and under local alternatives.

## 1. Introduction

The analysis of the unconditional distribution of linear predictors after model selection given in this paper complements and completes the results of Leeb [1], where the corresponding conditional distribution is considered, conditional on the outcome of the model selection step. The present paper builds on Leeb [1] as far as finite-sample results are concerned. For a large-sample analysis, however, we can not rely on that paper; the limit behavior of the unconditional cdf differs from that of the conditional cdfs so that a separate analysis is necessary. For a review of the relevant related literature and for an outline of applications of our results, we refer the reader to Leeb [1].

We consider a linear regression model $Y = X\theta + u$ with normal errors. (The normal linear model facilitates a detailed finite-sample analysis. Also note that asymptotic properties of the Gaussian location model can be generalized to a much larger model class including nonlinear models and models for dependent data, as long as appropriate standard regularity conditions guaranteeing asymptotic normality of the maximum likelihood estimator are satisfied.) We consider model selection by a sequence of 'general-to-specific' hypothesis tests; that is, starting from the overall model, a sequence of tests is used to simplify the model. The cdf of a linear function of the post-model-selection estimator (properly scaled and centered) is denoted by $G_{n,\theta,\sigma}(t)$. The notation suggests that this cdf depends on the sample size $n$, the regression parameter $\theta$, and on the error variance $\sigma^2$. An explicit formula for $G_{n,\theta,\sigma}(t)$ is given in (3.10) below. From this formula, we see that the distribution of, say, a linear predictor after model selection is significantly different from

[1]Department of Statistics, Yale University, 24 Hillhouse Avenue, New Haven, CT 06511.
* Research supported by the Max Kade Foundation and by the Austrian Science Foundation (FWF), project no. P13868-MAT. A preliminary version of this manuscript was written in February 2002.

(and more complex than) the Gaussian distribution that one would get without model selection. Because the cdf $G_{n,\theta,\sigma}(t)$ is quite difficult to analyze directly, we also provide a uniform asymptotic approximation to this cdf. This approximation, which we shall denote by $G^*_{n,\theta,\sigma}(t)$, is obtained by considering an 'idealized' scenario where the error variance $\sigma^2$ is treated as known and is used by the model selection procedure. The approximating cdf $G^*_{n,\theta,\sigma}(t)$ is much simpler and allows us to observe the main effects of model selection. Moreover, this approximation allows us to study the large-sample limit behavior of $G_{n,\theta,\sigma}(t)$ not only in the fixed-parameter case but also along sequences of parameters. The consideration of asymptotics along sequences of parameters is necessitated by a complication that seems to be inherent to post-model-selection estimators: Convergence of the finite-sample distributions to the large-sample limit distribution is non-uniform in the underlying parameters. (See Corollary 5.5 in Leeb and Pötscher [3], Appendix B in Leeb and Pötscher [4].) For applications like the computation of large-sample limit minimal coverage probabilities, it therefore appears to be necessary to study the limit behavior of $G_{n,\theta,\sigma}(t)$ along *sequences of parameters* $\theta^{(n)}$ and $\sigma^{(n)}$. We characterize all accumulation points of $G_{n,\theta^{(n)},\sigma^{(n)}}(t)$ for such sequences (with respect to weak convergence). Ex post, it turns out that, as far as possible accumulation points are concerned, it suffices to consider only a particular class of parameter sequences, namely local alternatives. Of course, the large-sample limit behavior of $G_{n,\theta,\sigma}(t)$ in the fixed-parameter case is contained in this analysis. Besides, we also consider the model selection probabilities, i.e., the probabilities of selecting each candidate model under consideration, in the finite-sample and in the large-sample limit case.

The remainder of the paper is organized as follows: In Section 2, we describe the basic framework of our analysis and the quantities of interest: The post-model-selection estimator $\tilde{\theta}$ and the cdf $G_{n,\theta,\sigma}(t)$. Besides, we also introduce the 'idealized post-model-selection estimator' $\tilde{\theta}^*$ and the cdf $G^*_{n,\theta,\sigma}(t)$, which correspond to the case where the error variance is known. In Section 3, we derive finite-sample expansions of the aforementioned cdfs, and we discuss and illustrate the effects of the model selection step in finite samples. Section 4 contains an approximation result which shows that $G_{n,\theta,\sigma}(t)$ and $G^*_{n,\theta,\sigma}(t)$ are asymptotically uniformly close to each other. With this, we can analyze the large-sample limit behavior of the two cdfs in Section 5. All proofs are relegated to the appendices.

## 2. The model and estimators

Consider the linear regression model

$$(2.1) \qquad Y = X\theta + u,$$

where $X$ is a non-stochastic $n \times P$ matrix with $\text{rank}(X) = P$ and $u \sim N(0, \sigma^2 I_n)$, $\sigma^2 > 0$. Here $n$ denotes the sample size and we assume $n > P \geq 1$. In addition, we assume that $Q = \lim_{n\to\infty} X'X/n$ exists and is non-singular (this assumption is not needed in its full strength for all of the asymptotic results; cf. Remark 2.1). Similarly as in Pötscher [6], we consider model selection from a collection of nested models $M_{\mathcal{O}} \subseteq M_{\mathcal{O}+1} \subseteq \cdots \subseteq M_P$ which are given by $M_p = \{(\theta_1, \ldots, \theta_P)' \in \mathbb{R}^P : \theta_{p+1} = \cdots = \theta_P = 0\}$ $(0 \leq p \leq P)$. Hence, the model $M_p$ corresponds to the situation where only the first $p$ regressors in (2.1) are included. For the most parsimonious model under consideration, i.e., for $M_{\mathcal{O}}$, we assume that $\mathcal{O}$ satisfies $0 \leq \mathcal{O} < P$; if $\mathcal{O} > 0$, this model contains those components of the parameter that will not be subject to model selection. Note that $M_0 = \{(0, \ldots, 0)'\}$ and $M_P = \mathbb{R}^P$. We call $M_p$ the regression model of order $p$.

The following notation will prove useful. For matrices $B$ and $C$ of the same row-dimension, the column-wise concatenation of $B$ and $C$ is denoted by $(B:C)$. If $D$ is an $m \times P$ matrix, let $D[p]$ denote the matrix of the first $p$ columns of $D$. Similarly, let $D[\neg p]$ denote the matrix of the last $P - p$ columns of $D$. If $x$ is a $P \times 1$ (column-) vector, we write in abuse of notation $x[p]$ and $x[\neg p]$ for $(x'[p])'$ and $(x'[\neg p])'$, respectively. (We shall use these definitions also in the 'boundary' cases $p = 0$ and $p = P$. It will always be clear from the context how expressions like $D[0]$, $D[\neg P]$, $x[0]$, or $x[\neg P]$ are to be interpreted.) As usual the $i$-th component of a vector $x$ will be denoted by $x_i$; in a similar fashion, denote the entry in the $i$-th row and $j$-th column of a matrix $B$ by $B_{i,j}$.

The restricted least-squares estimator for $\theta$ under the restriction $\theta[\neg p] = 0$ will be denoted by $\tilde{\theta}(p)$, $0 \leq p \leq P$ (in case $p = P$, the restriction is void, of course). Note that $\tilde{\theta}(p)$ is given by the $P \times 1$ vector whose first $p$ components are given by $(X[p]'X[p])^{-1}X[p]'Y$, and whose last $P - p$ components are equal to zero; the expressions $\tilde{\theta}(0)$ and $\tilde{\theta}(P)$, respectively, are to be interpreted as the zero-vector in $\mathbb{R}^P$ and as the unrestricted least-squares estimator for $\theta$. Given a parameter vector $\theta$ in $\mathbb{R}^P$, the order of $\theta$, relative to the set of models $M_0, \ldots, M_P$, is defined as $p_0(\theta) = \min\{p: 0 \leq p \leq P, \theta \in M_p\}$. Hence, if $\theta$ is the true parameter vector, only models $M_p$ of order $p \geq p_0(\theta)$ are correct models, and $M_{p_0(\theta)}$ is the most parsimonious correct model for $\theta$ among $M_0, \ldots, M_P$. We stress that $p_0(\theta)$ is a property of a *single parameter*, and hence needs to be distinguished from the notion of the order of the model $M_p$ introduced earlier, which is a property of the *set of parameters $M_p$*.

A model selection procedure in general is now nothing else than a data-driven (measurable) rule $\hat{p}$ that selects a value from $\{\mathcal{O}, \ldots, P\}$ and thus selects a model from the list of candidate models $M_\mathcal{O}, \ldots, M_P$. In this paper, we shall consider a model selection procedure based on a sequence of 'general-to-specific' hypothesis tests, which is given as follows: The sequence of hypotheses $H_0^p : p_0(\theta) < p$ is tested against the alternatives $H_1^p : p_0(\theta) = p$ in decreasing order starting at $p = P$. If, for some $p > \mathcal{O}$, $H_0^p$ is the first hypothesis in the process that is rejected, we set $\hat{p} = p$. If no rejection occurs until even $H_0^{\mathcal{O}+1}$ is accepted, we set $\hat{p} = \mathcal{O}$. Each hypothesis in this sequence is tested by a kind of $t$-test where the error variance is always estimated from the overall model. More formally, we have

$$\hat{p} = \max\{p: |T_p| \geq c_p, 0 \leq p \leq P\},$$

where the test-statistics are given by $T_0 = 0$ and by $T_p = \sqrt{n}\tilde{\theta}_p(p)/(\hat{\sigma}\xi_{n,p})$ with

$$(2.2) \qquad \xi_{n,p} = \left(\left[\left(\frac{X[p]'X[p]}{n}\right)^{-1}\right]_{p,p}\right)^{\frac{1}{2}} \qquad (0 < p \leq P)$$

being the (non-negative) square root of the $p$-th diagonal element of the matrix indicated, and with $\hat{\sigma}^2 = (n-P)^{-1}(Y - X\tilde{\theta}(P))'(Y - X\tilde{\theta}(P))$ (cf. also Remark 6.2 in Leeb [1] concerning other variance estimators). The critical values $c_p$ are independent of sample size (cf., however, Remark 2.1) and satisfy $0 < c_p < \infty$ for $\mathcal{O} < p \leq P$. We also set $c_\mathcal{O} = 0$ in order to restrict $\hat{p}$ to the range of candidate models under consideration, i.e., to $\{\mathcal{O}, \mathcal{O}+1, \ldots, P\}$. Note that under the hypothesis $H_0^p$ the statistic $T_p$ is $t$-distributed with $n - P$ degrees of freedom for $0 < p \leq P$. The so defined model selection procedure $\hat{p}$ is conservative (or over-consistent): The probability of selecting an incorrect model, i.e., the probability of the event $\{\hat{p} < p_0(\theta)\}$, converges to zero as the sample size increases; the probability of selecting a correct (but possibly over-parameterized) model, i.e., the probability of the event $\{\hat{p} = p\}$ for $p$ satisfying $\max\{p_0(\theta), \mathcal{O}\} \leq p \leq P$, converges to a positive limit; cf. (5.6) below.

The post-model-selection estimator $\tilde{\theta}$ is now defined as follows: On the event $\hat{p} = p$, $\tilde{\theta}$ is given by the restricted least-squares estimator $\tilde{\theta}(p)$, i.e.,

$$(2.3) \qquad \tilde{\theta} = \sum_{p=\mathcal{O}}^{P} \tilde{\theta}(p) \mathbf{1}\{\hat{p} = p\}.$$

To study the distribution of a linear function of $\tilde{\theta}$, let $A$ be a non-stochastic $k \times P$ matrix of rank $k$ ($1 \leq k \leq P$). Examples for $A$ include the case where $A$ equals a $1 \times P$ (row-) vector $x_f$ if the object of interest is the linear predictor $x_f \tilde{\theta}$, or the case where $A = (I_s : 0)$, say, if the object of interest is an $s \times 1$ subvector of $\theta$. We shall consider the cdf

$$(2.4) \qquad G_{n,\theta,\sigma}(t) = P_{n,\theta,\sigma}\left(\sqrt{n} A(\tilde{\theta} - \theta) \leq t\right) \qquad (t \in \mathbb{R}^k).$$

Here and in the following, $P_{n,\theta,\sigma}(\cdot)$ denotes the probability measure corresponding to a sample of size $n$ from (2.1) under the true parameters $\theta$ and $\sigma$. For convenience we shall refer to (2.4) as the cdf of $A\tilde{\theta}$, although (2.4) is in fact the cdf of an affine transformation of $A\tilde{\theta}$.

For theoretical reasons we shall also be interested in the idealized model selection procedure which assumes knowledge of $\sigma^2$ and hence uses $T_p^*$ instead of $T_p$, where $T_p^* = \sqrt{n}\tilde{\theta}_p(p)/(\sigma \xi_{n,p})$, $0 < p \leq P$, and $T_0^* = 0$. The corresponding model selector is denoted by $\hat{p}^*$ and the resulting idealized 'post-model-selection estimator' by $\tilde{\theta}^*$. Note that under the hypothesis $H_0^p$ the variable $T_p^*$ is standard normally distributed for $0 < p \leq P$. The corresponding cdf will be denoted by $G_{n,\theta,\sigma}^*(t)$, i.e.,

$$(2.5) \qquad G_{n,\theta,\sigma}^*(t) = P_{n,\theta,\sigma}\left(\sqrt{n} A(\tilde{\theta}^* - \theta) \leq t\right) \qquad (t \in \mathbb{R}^k).$$

For convenience we shall also refer to (2.5) as the cdf of $A\tilde{\theta}^*$.

**Remark 2.1.** Some of the assumptions introduced above are made only to simplify the exposition and can hence easily be relaxed. This includes, in particular, the assumption that the critical values $c_p$ used by the model selection procedure do not depend on sample size, and the assumption that the regressor matrix $X$ is such that $X'X/n$ converges to a positive definite limit $Q$ as $n \to \infty$. For the finite-sample results in Section 3 below, these assumptions are clearly inconsequential. Moreover, for the large-sample limit results in Sections 4 and 5 below, these assumptions can be relaxed considerably. For the details, see Remark 6.1(i)–(iii) in Leeb [1], which also applies, mutatis mutandis, to the results in the present paper.

## 3. Finite-sample results

Some further preliminaries are required before we can proceed. The expected value of the restricted least-squares estimator $\tilde{\theta}(p)$ will be denoted by $\eta_n(p)$ and is given by the $P \times 1$ vector

$$(3.1) \quad \eta_n(p) = \begin{pmatrix} \theta[p] + (X[p]'X[p])^{-1}X[p]'X[\neg p]\theta[\neg p] \\ (0,\ldots,0)' \end{pmatrix}$$

with the conventions that $\eta_n(0) = (0,\ldots,0)' \in \mathbb{R}^P$ and $\eta_n(P) = \theta$. Furthermore, let $\Phi_{n,p}(t)$, $t \in \mathbb{R}^k$, denote the cdf of $\sqrt{n}A(\tilde{\theta}(p) - \eta_n(p))$, i.e., $\Phi_{n,p}(t)$ is the cdf of a centered Gaussian random vector with covariance matrix $\sigma^2 A[p](X[p]'X[p]/n)^{-1}A[p]'$ in case $p > 0$, and the cdf of point-mass at zero in $\mathbb{R}^k$ in case $p = 0$. If $p > 0$ and if the matrix $A[p]$ has rank $k$, then $\Phi_{n,p}(t)$ has a density with respect to Lebesgue measure, and we shall denote this density by $\phi_{n,p}(t)$. We note that $\eta_n(p)$ depends on $\theta$ and that $\Phi_{n,p}(t)$ depends on $\sigma$ (in case $p > 0$), although these dependencies are not shown explicitly in the notation.

For $p > 0$, the conditional distribution of $\sqrt{n}\tilde{\theta}_p(p)$ given $\sqrt{n}A(\tilde{\theta}(p) - \eta_n(p)) = z$ is a Gaussian distribution with mean $\sqrt{n}\eta_{n,p}(p) + b_{n,p}z$ and variance $\sigma^2 \zeta_{n,p}^2$, where

$$(3.2) \quad b_{n,p} = C_n^{(p)'}(A[p](X[p]'X[p]/n)^{-1}A[p]')^{-}, \text{ and}$$

$$(3.3) \quad \zeta_{n,p}^2 = \xi_{n,p}^2 - b_{n,p}C_n^{(p)}.$$

In the displays above, $C_n^{(p)}$ stands for $A[p](X[p]'X[p]/n)^{-1}e_p$, with $e_p$ denoting the $p$-th standard basis vector in $\mathbb{R}^p$, and $(A[p](X[p]'X[p]/n)^{-1}A[p]')^{-}$ denotes a generalized inverse of the matrix indicated (cf. Note 3(v) in Section 8a.2 of Rao [7]). Note that, in general, the quantity $b_{n,p}z$ depends on the choice of generalized inverse in (3.2); however, for $z$ in the column-space of $A[p]$, $b_{n,p}z$ is invariant under the choice of inverse; cf. Lemma A.2 in Leeb [1]. Since $\sqrt{n}A(\tilde{\theta}(p) - \eta_n(p))$ lies in the column-space of $A[p]$, the conditional distribution of $\sqrt{n}\tilde{\theta}_p(p)$ given $\sqrt{n}A(\tilde{\theta}(p) - \eta_n(p)) = z$ is thus well-defined by the above. Observe that the vector of covariances between $A\tilde{\theta}(p)$ and $\tilde{\theta}_p(p)$ is given by $\sigma^2 n^{-1} C_n^{(p)}$. In particular, note that $A\tilde{\theta}(p)$ and $\tilde{\theta}_p(p)$ are uncorrelated if and only if $\zeta_{n,p}^2 = \xi_{n,p}^2$ (or, equivalently, if and only if $b_{n,p}z = 0$ for all $z$ in the column-space of $A[p]$); again, see Lemma A.2 in Leeb [1].

Finally, for $M$ denoting a univariate Gaussian random variable with zero mean and variance $s^2 \geq 0$, we abbreviate the probability $P(|M - a| < b)$ by $\Delta_s(a,b)$, $a \in \mathbb{R} \cup \{-\infty, \infty\}$, $b \in \mathbb{R}$. Note that $\Delta_s(\cdot, \cdot)$ is symmetric around zero in its first argument, and that $\Delta_s(-\infty, b) = \Delta_s(\infty, b) = 0$ holds. In case $s = 0$, $M$ is to be interpreted as being equal to zero, such that $\Delta_0(a,b)$ equals one if $|a| < b$ and zero otherwise; i.e., $\Delta_0(a,b)$ reduces to an indicator function.

### 3.1. The known-variance case

The cdf $G^*_{n,\theta,\sigma}(t)$ can be expanded as a weighted sum of conditional cdfs, conditional on the outcome of the model selection step, where the weights are given by the corresponding model selection probabilities. To this end, let $G^*_{n,\theta,\sigma}(t|p)$ denote the conditional cdf of $\sqrt{n}A(\tilde{\theta}^* - \theta)$ given that $\hat{p}^*$ equals $p$ for $\mathcal{O} \leq p \leq P$; that is, $G^*_{n,\theta,\sigma}(t|p) = P_{n,\theta,\sigma}(\sqrt{n}A(\tilde{\theta}^* - \theta) \leq t \mid \hat{p}^* = p)$, with $t \in \mathbb{R}^k$. Moreover, let $\pi^*_{n,\theta,\sigma}(p) = P_{n,\theta,\sigma}(\hat{p}^* = p)$ denote the corresponding model selection probability. Then the unconditional cdf $G^*_{n,\theta,\sigma}(t)$ can be written as

$$(3.4) \quad G^*_{n,\theta,\sigma}(t) = \sum_{p=\mathcal{O}}^{P} G^*_{n,\theta,\sigma}(t|p)\pi^*_{n,\theta,\sigma}(p).$$

Explicit finite-sample formulas for $G^*_{n,\theta,\sigma}(t|p)$, $\mathcal{O} \leq p \leq P$, are given in Leeb [1], equations (10) and (13). Let $\gamma(\xi_{n,q}, s) = \Delta_{\sigma\xi_{n,q}}(\sqrt{n}\eta_{n,q}(q), sc_q\sigma\xi_{n,q})$, and $\gamma^*(\zeta_{n,q}, z, s) = \Delta_{\sigma\zeta_{n,q}}(\sqrt{n}\eta_{n,q}(q) + b_{n,q}z, sc_q\sigma\xi_{n,q})$ It is elementary to verify that $\pi^*_{n,\theta,\sigma}(\mathcal{O})$ is given by

$$(3.5) \qquad \pi^*_{n,\theta,\sigma}(\mathcal{O}) = \prod_{q=\mathcal{O}+1}^{P} \gamma(\xi_{n,q}, 1)$$

while, for $p > \mathcal{O}$, we have

$$(3.6) \qquad \pi^*_{n,\theta,\sigma}(p) = (1 - \gamma(\xi_{n,p}, 1)) \times \prod_{q=p+1}^{P} \gamma(\xi_{n,q}, 1).$$

(This follows by arguing as in the discussion leading up to (12) of Leeb [1], and by using Proposition 3.1 of that paper.) Observe that the model selection probability $\pi^*_{n,\theta,\sigma}(p)$ is always positive for each $p$, $\mathcal{O} \leq p \leq P$.

Plugging the formulas for the conditional cdfs obtained in Leeb [1] and the above formulas for the model selection probabilities into (3.4), we obtain that $G^*_{n,\theta,\sigma}(t)$ is given by

$$(3.7) \qquad \begin{aligned} G^*_{n,\theta,\sigma}(t) &= \Phi_{n,\mathcal{O}}(t - \sqrt{n}A(\eta_n(\mathcal{O}) - \theta)) \prod_{q=\mathcal{O}+1}^{P} \gamma(\xi_{n,q}, 1) \\ &\quad + \sum_{p=\mathcal{O}+1}^{P} \int_{z \leq t - \sqrt{n}A(\eta_n(p) - \theta)} (1 - \gamma^*(\zeta_{n,p}, z, 1)) \Phi_{n,p}(dz) \\ &\quad \times \prod_{q=p+1}^{P} \gamma(\xi_{n,q}, 1). \end{aligned}$$

In the above display, $\Phi_{n,p}(dz)$ denotes integration with respect to the measure induced by the cdf $\Phi_{n,p}(t)$ on $\mathbb{R}^k$.

### 3.2. The unknown-variance case

Similar to the known-variance case, define $G_{n,\theta,\sigma}(t|p) = P_{n,\theta,\sigma}(\sqrt{n}A(\tilde{\theta} - \theta) \leq t | \hat{p} = p)$ and $\pi_{n,\theta,\sigma}(p) = P_{n,\theta,\sigma}(\hat{p} = p)$, $\mathcal{O} \leq p \leq P$. Then $G_{n,\theta,\sigma}(t)$ can be expanded as the sum of the terms $G_{n,\theta,\sigma}(t|p)\pi_{n,\theta,\sigma}(p)$ for $p = \mathcal{O}, \ldots, P$, similar to (3.4).

For the model selection probabilities, we argue as in Section 3.2 of Leeb and Pötscher [3] to obtain that $\pi_{n,\theta,\sigma}(\mathcal{O})$ equals

$$(3.8) \qquad \pi_{n,\theta,\sigma}(\mathcal{O}) = \int_0^\infty \prod_{q=\mathcal{O}+1}^{P} \gamma(\xi_{n,q}, s) h(s) ds,$$

where $h$ denotes the density of $\hat{\sigma}/\sigma$, i.e., $h$ is the density of $(n - P)^{-1/2}$ times the square-root of a chi-square distributed random variable with $n - P$ degrees of freedom. In a similar fashion, for $p > \mathcal{O}$, we get

$$(3.9) \qquad \pi_{n,\theta,\sigma}(p) = \int_0^\infty (1 - \gamma(\xi_{n,p}, s)) \prod_{q=p+1}^{P} \gamma(\xi_{n,q}, s) h(s) ds;$$

cf. the argument leading up to (18) in Leeb [1]. As in the known-variance case, the model selection probabilities are all positive.

Using the formulas for the conditional cdfs $G_{n,\theta,\sigma}(t|p)$, $\mathcal{O} \leq p \leq P$, given in Leeb [1], equations (14) and (16)–(18), the unconditional cdf $G_{n,\theta,\sigma}(t)$ is thus seen to be given by

$$
\begin{aligned}
G_{n,\theta,\sigma}(t) &= \Phi_{n,\mathcal{O}}(t - \sqrt{n}A(\eta_n(\mathcal{O}) - \theta)) \int_0^\infty \prod_{q=\mathcal{O}+1}^P \gamma(\xi_{n,q}, s) h(s) ds \\
&\quad + \sum_{p=\mathcal{O}+1}^P \int_{z \leq t - \sqrt{n}A(\eta_n(p) - \theta)} \left[ \int_0^\infty (1 - \gamma^*(\zeta_{n,p}, z, s)) \right. \\
&\quad \left. \times \prod_{q=p+1}^P \gamma(\xi_{n,q}, s) h(s) ds \right] \Phi_{n,p}(dz).
\end{aligned}
$$
(3.10)

Observe that $G_{n,\theta,\sigma}(t)$ is in fact a smoothed version of $G^*_{n,\theta,\sigma}(t)$: Indeed, the right-hand side of formula (3.10) for $G_{n,\theta,\sigma}(t)$ is obtained by taking the right-hand side of formula (3.7) for $G^*_{n,\theta,\sigma}(t)$, changing the last argument of $\gamma(\xi_{n,q}, 1)$ and $\gamma^*(\zeta_{n,q}, z, 1)$ from 1 to $s$ for $q = \mathcal{O} + 1, \ldots, P$, integrating with respect to $h(s)ds$, and interchanging the order of integration. Similar considerations apply, mutatis mutandis, to the model selection probabilities $\pi_{n,\theta,\sigma}(p)$ and $\pi^*_{n,\theta,\sigma}(p)$ for $\mathcal{O} \leq p \leq P$.

### 3.3. Discussion

#### 3.3.1. General Observations

The cdfs $G^*_{n,\theta,\sigma}(t)$ and $G_{n,\theta,\sigma}(t)$ need not have densities with respect to Lebesgue measure on $\mathbb{R}^k$. However, densities do exist if $\mathcal{O} > 0$ and the matrix $A[\mathcal{O}]$ has rank $k$. In that case, the density of $G_{n,\theta,\sigma}(t)$ is given by

$$
\begin{aligned}
&\phi_{n,\mathcal{O}}(t - \sqrt{n}A(\eta_n(\mathcal{O}) - \theta)) \int_0^\infty \prod_{q=\mathcal{O}+1}^P \gamma(\xi_{n,q}, s) h(s) ds \\
&\quad + \sum_{p=\mathcal{O}+1}^P \left[ \int_0^\infty (1 - \gamma^*(\zeta_{n,p}, t - \sqrt{n}A(\eta_n(p) - \theta), s)) \right. \\
&\quad \left. \times \prod_{q=p+1}^P \gamma(\xi_{n,q}, s) h(s) ds \right] \phi_{n,p}(t - \sqrt{n}A(\eta_n(p) - \theta)).
\end{aligned}
$$
(3.11)

(Given that $\mathcal{O} > 0$ and that $A[\mathcal{O}]$ has rank $k$, we see that $A[p]$ has rank $k$ and that the Lebesgue density $\phi_{n,p}(t)$ of $\Phi_{n,p}(t)$ exists for each $p = \mathcal{O}, \ldots, P$. We hence may write the integrals with respect to $\Phi_{n,p}(dz)$ in (3.10) as integrals with respect to $\phi_{n,p}(z)dz$. Differentiating the resulting formula for $G_{n,\theta,\sigma}(t)$ with respect to $t$, we get (3.11).) Similarly, the Lebesgue density of $G^*_{n,\theta,\sigma}(t)$ can be obtained by differentiating the right-hand side of (3.7), provided that $\mathcal{O} > 0$ and $A[\mathcal{O}]$ has rank $k$. Conversely, if that condition is violated, then some of the conditional cdfs are degenerate and Lebesgue densities do not exist. (Note that on the event $\hat{p} = p$, $A\tilde{\theta}$ equals $A\tilde{\theta}(p)$, and recall that the last $P - p$ coordinates of $\tilde{\theta}(p)$ are constant equal to zero. Therefore $A\tilde{\theta}(0)$ is the zero-vector in $\mathbb{R}^k$ and, for $p > 0$, $A\tilde{\theta}(p)$ is concentrated in the column space of $A[p]$. On the event $\hat{p}^* = p$, a similar argument applies to $A\tilde{\theta}^*$.)

Both cdfs $G_{n,\theta,\sigma}^*(t)$ and $G_{n,\theta,\sigma}(t)$ are given by a weighted sum of conditional cdfs, cf. (3.7) and (3.10), where the weights are given by the model-selection probabilities (which are always positive in finite samples). For a detailed discussion of the conditional cdfs, the reader is referred to Section 3.3 of Leeb [1].

The cdf $G_{n,\theta,\sigma}(t)$ is typically highly non-Gaussian. A notable exception where $G_{n,\theta,\sigma}(t)$ reduces to the Gaussian cdf $\Phi_{n,P}(t)$ for each $\theta \in \mathbb{R}^P$ occurs in the special case where $\tilde{\theta}_p(p)$ is uncorrelated with $A\tilde{\theta}(p)$ for each $p = \mathcal{O}+1, \ldots, P$. In this case, we have $A\tilde{\theta}(p) = A\tilde{\theta}(P)$ for each $p = \mathcal{O}, \ldots, P$ (cf. the discussion following (20) in Leeb [1]). From this and in view of (2.3), it immediately follows that $G_{n,\theta,\sigma}(t) = \Phi_{n,P}(t)$, independent of $\theta$ and $\sigma$. (The same considerations apply, mutatis mutandis, to $G_{n,\theta,\sigma}^*(t)$.) Clearly, this case is rather special, because it entails that fitting the overall model with $P$ regressors gives the same estimator for $A\theta$ as fitting the restricted model with $\mathcal{O}$ regressors only.

To compare the distribution of a linear function of the post-model-selection estimator with the distribution of the post-model-selection estimator itself, note that the cdf of $\tilde{\theta}$ can be studied in our framework by setting $A$ equal to $I_P$ (and $k$ equal to $P$). Obviously, the distribution of $\tilde{\theta}$ does not have a density with respect to Lebesgue measure. Moreover, $\tilde{\theta}_p(p)$ is always perfectly correlated with $\tilde{\theta}(p)$ for each $p = 1, \ldots, P$, such that the special case discussed above can not occur (for $A$ equal to $I_P$).

### 3.3.2. An illustrative example

We now exemplify the possible shapes of the finite-sample distributions in a simple setting. To this end, we set $P = 2$, $\mathcal{O} = 1$, $A = (1,0)$, and $k = 1$ for the rest of this section. The choice of $P = 2$ gives a special case of the model (2.1), namely

$$(3.12) \qquad Y_i = \theta_1 X_{i,1} + \theta_2 X_{i,2} + u_i \qquad (1 \le i \le n).$$

With $\mathcal{O} = 1$, the first regressor is always included in the model, and a pre-test will be employed to decide whether or not to include the second one. The two model selectors $\hat{p}$ and $\hat{p}^*$ thus decide between two candidate models, $M_1 = \{(\theta_1, \theta_2)' \in \mathbb{R}^2 : \theta_2 = 0\}$ and $M_2 = \{(\theta_1, \theta_2)' \in \mathbb{R}^2\}$. The critical value for the test between $M_1$ and $M_2$, i.e., $c_2$, will be chosen later (recall that we have set $c_\mathcal{O} = c_1 = 0$). With our choice of $A = (1,0)$, we see that $G_{n,\theta,\sigma}(t)$ and $G_{n,\theta,\sigma}^*(t)$ are the cdfs of $\sqrt{n}(\tilde{\theta}_1 - \theta_1)$ and $\sqrt{n}(\tilde{\theta}_1^* - \theta_1)$, respectively.

Since the matrix $A[\mathcal{O}]$ has rank one and $k = 1$, the cdfs of $\sqrt{n}(\tilde{\theta}_1 - \theta_1)$ and $\sqrt{n}(\tilde{\theta}_1^* - \theta_1)$ both have Lebesgue densities. To obtain a convenient expression for these densities, we write $\sigma^2(X'X/n)^{-1}$, i.e., the covariance matrix of the least-squares estimator based on the overall model (3.12), as

$$\sigma^2 \left(\frac{X'X}{n}\right)^{-1} = \begin{pmatrix} \sigma_1^2 & \sigma_{1,2} \\ \sigma_{1,2} & \sigma_2^2 \end{pmatrix}.$$

The elements of this matrix depend on sample size $n$, but we shall suppress this dependence in the notation. It will prove useful to define $\rho = \sigma_{1,2}/(\sigma_1 \sigma_2)$, i.e., $\rho$ is the correlation coefficient between the least-squares estimators for $\theta_1$ and $\theta_2$ in model (3.12). Note that here we have $\phi_{n,2}(t) = \sigma_1^{-1} \phi(t/\sigma_1)$ and $\phi_{n,1}(t) =$ $\sigma_1^{-1}(1-\rho^2)^{-1/2}\phi(t(1-\rho^2)^{-1/2}/\sigma_1)$ with $\phi(t)$ denoting the univariate standard Gaussian density. The density of $\sqrt{n}(\tilde{\theta}_1 - \theta_1)$ is given by

$$
(3.13) \quad \begin{aligned} & \phi_{n,1}(t + \sqrt{n}\theta_2\rho\sigma_1/\sigma_2)\int_0^\infty \Delta_1(\sqrt{n}\theta_2/\sigma_2, sc_2)h(s)ds \\ & + \phi_{n,2}(t)\int_0^\infty (1 - \Delta_1(\frac{\sqrt{n}\theta_2/\sigma_2 + \rho t/\sigma_1}{\sqrt{1-\rho^2}}, \frac{sc_2}{\sqrt{1-\rho^2}}))h(s)ds; \end{aligned}
$$

recall that $\Delta_1(a,b)$ is equal to $\Phi(a+b) - \Phi(a-b)$, where $\Phi(t)$ denotes the standard univariate Gaussian cdf, and note that here $h(s)$ denotes the density of $(n-2)^{-1/2}$ times the square-root of a chi-square distributed random variable with $n-2$ degrees of freedom. Similarly, the density of $\sqrt{n}(\tilde{\theta}_1^* - \theta_1)$ is given by

$$
(3.14) \quad \begin{aligned} & \phi_{n,1}(t + \sqrt{n}\theta_2\rho\sigma_1/\sigma_2)\Delta_1(\sqrt{n}\theta_2/\sigma_2, c_2) \\ & + \phi_{n,2}(t)(1 - \Delta_1(\frac{\sqrt{n}\theta_2/\sigma_2 + \rho t/\sigma_1}{\sqrt{1-\rho^2}}, \frac{c_2}{\sqrt{1-\rho^2}})). \end{aligned}
$$

Note that both densities depend on the regression parameter $(\theta_1, \theta_2)'$ only through $\theta_2$, and that these densities depend on the error variance $\sigma^2$ and on the regressor matrix $X$ only through $\sigma_1$, $\sigma_2$, and $\rho$. Also note that the expressions in (3.13) and (3.14) are unchanged if $\rho$ is replaced by $-\rho$ and, at the same time, the argument $t$ is replaced by $-t$. Similarly, replacing $\theta_2$ and $t$ by $-\theta_2$ and $-t$, respectively, leaves (3.13) and (3.14) unchanged. The same applies also to the conditional densities considered below; cf. (3.15) and (3.16). We therefore consider only non-negative values of $\rho$ and $\theta_2$ in the numerical examples below.

From (3.14) we can also read-off the conditional densities of $\sqrt{n}(\tilde{\theta}_1^* - \theta_1)$, conditional on selecting the model $M_p$ for $p = 1$ and $p = 2$, which will be useful later: The unconditional cdf of $\sqrt{n}(\tilde{\theta}_1^* - \theta_1)$ is the weighted sum of two conditional cdfs, conditional on selecting the model $M_1$ and $M_2$, respectively, weighted by the corresponding model selection probabilities; cf. (3.4) and the attending discussion. Hence, the unconditional density is the sum of the conditional densities multiplied by the corresponding model selection probabilities. In the simple setting considered here, the probability of $\hat{p}^*$ selecting $M_1$, i.e., $\pi^*_{n,\theta,\sigma}(1)$, equals $\Delta_1(\sqrt{n}\theta_2/\sigma_2, c_2)$ in view of (3.5) and because $\mathcal{O} = 1$, and $\pi^*_{n,\theta,\sigma}(2) = 1 - \pi^*_{n,\theta,\sigma}(1)$. Thus, conditional on selecting the model $M_1$, the density of $\sqrt{n}(\tilde{\theta}_1^* - \theta_1)$ is given by

$$(3.15) \quad \phi_{n,1}(t + \sqrt{n}\theta_2\rho\sigma_1/\sigma_2).$$

Conditional on selecting $M_2$, the density of $\sqrt{n}(\tilde{\theta}_1^* - \theta_1)$ equals

$$(3.16) \quad \phi_{n,2}(t)\frac{1 - \Delta_1((\sqrt{n}\theta_2/\sigma_2 + \rho t/\sigma_1)/\sqrt{1-\rho^2}, c_2/\sqrt{1-\rho^2})}{1 - \Delta_1(\sqrt{n}\theta_2/\sigma_2, c_2)}.$$

This can be viewed as a 'deformed' version of $\phi_{n,2}(t)$, i.e., the density of $\sqrt{n}(\tilde{\theta}_1(2) - \theta_1)$, where the deformation is governed by the fraction in (3.16). The conditional densities of $\sqrt{n}(\tilde{\theta}_1 - \theta_1)$ can be obtained and interpreted in a similar fashion from (3.13), upon observing that $\pi_{n,\theta,\sigma}(1)$ here equals $\int_0^\infty \Delta_1(\sqrt{n}\theta_2/\sigma_2, sc_2)h(s)ds$ in view of (3.8)

Figure 1 illustrates some typical shapes of the densities of $\sqrt{n}(\tilde{\theta}_1 - \theta_1)$ and $\sqrt{n}(\tilde{\theta}_1^* - \theta_1)$ given in (3.13) and (3.14), respectively, for $\rho = 0.75$, $n = 7$, and for various values of $\theta_2$. Note that the densities of $\sqrt{n}(\tilde{\theta}_1 - \theta_1)$ and $\sqrt{n}(\tilde{\theta}_1^* - \theta_1)$, corresponding to the unknown-variance case and the (idealized) known-variance case, are very close to each other. In fact, the small sample size, i.e., $n = 7$, was chosen because for larger $n$ these two densities are visually indistinguishable in plots as in Figure 1 (this phenomenon is analyzed in detail in the next section). For $\theta_2 = 0$ in Figure 1, the density of $\sqrt{n}(\tilde{\theta}_1^* - \theta_1)$, although seemingly close to being Gaussian, is in fact a mixture of a Gaussian density and a bimodal density; this is explained in detail below. For the remaining values of $\theta_2$ considered in Figure 1, the density of $\sqrt{n}(\tilde{\theta}_1^* - \theta_1)$ is clearly non-Gaussian, namely skewed in case $\theta_2 = 0.1$, bimodal in case $\theta_2 = 0.75$, and highly non-symmetric in case $\theta_2 = 1.2$. Overall, we see that the finite-sample density of $\sqrt{n}(\tilde{\theta}_1^* - \theta_1)$ can exhibit a variety of different shapes. Exactly the same applies to the density of $\sqrt{n}(\tilde{\theta}_1 - \theta_1)$. As a point of interest, we note that these different shapes occur for values of $\theta_2$ in a quite narrow range: For example, in the setting of Figure 1, the uniformly most powerful test of the hypothesis $\theta_2 = 0$ against $\theta_2 > 0$ with level $0.95$, i.e., a one-sided t-test, has a power of only $0.27$ at the alternative $\theta_2 = 1.2$. This suggests that estimating the distribution of $\sqrt{n}(\tilde{\theta}_1 - \theta_1)$ is difficult here. (See also Leeb and Pöstcher [4] as well as Leeb and Pöstcher [2] for a thorough analysis of this difficulty.)

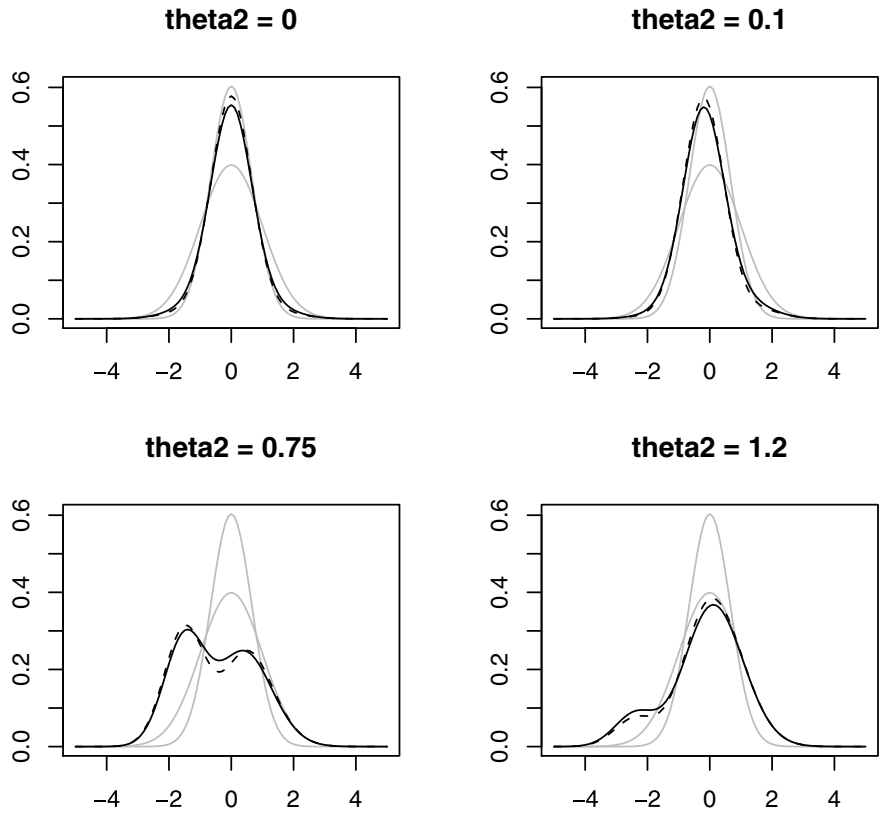

FIG 1. *The densities of $\sqrt{n}(\tilde{\theta}_1 - \theta_1)$ (black solid line) and of $\sqrt{n}(\tilde{\theta}_1^* - \theta_1)$ (black dashed line) for the indicated values of $\theta_2$, $n = 7$, $\rho = 0.75$, and $\sigma_1 = \sigma_2 = 1$. The critical value of the test between $M_1$ and $M_2$ was set to $c_2 = 2.015$, corresponding to a t-test with significance level $0.9$. For reference, the gray curves are Gaussian densities $\phi_{n,1}(t)$ (larger peak) and $\phi_{n,2}(t)$ (smaller peak).*

We stress that the phenomena shown in Figure 1 are not caused by the small sample size, i.e., $n = 7$. This becomes clear upon inspection of (3.13) and (3.14), which depend on $\theta_2$ through $\sqrt{n}\theta_2$ (for fixed $\sigma_1$, $\sigma_2$ and $\rho$). Hence, for other values of $n$, one obtains plots essentially similar to Figure 1, provided that the range of values of $\theta_2$ is adapted accordingly.

We now show how the shape of the unconditional densities can be explained by the shapes of the conditional densities together with the model selection probabilities. Since the unknown-variance case and the known-variance case are very similar as seen above, we focus on the latter. In Figure 2 below, we give the conditional densities of $\sqrt{n}(\tilde{\theta}_1^* - \theta_1)$, conditional on selecting the model $M_p$, $p = 1, 2$, cf. (3.15) and (3.16), and the corresponding model selection probabilities in the same setting as in Figure 1.

The unconditional densities of $\sqrt{n}(\tilde{\theta}_1^* - \theta_1)$ in each panel of Figure 1 are the sum of the two conditional densities in the corresponding panel in Figure 2, weighted by the model selection probabilities, i.e, $\pi_{n,\theta,\sigma}^*(1)$ and $\pi_{n,\theta,\sigma}^*(2)$. In other words, in each panel of Figure 2, the solid black curve gets the weight given in parentheses, and the dashed black curve gets one minus that weight. In case $\theta_2 = 0$, the probability of selecting model $M_1$ is very large, and the corresponding conditional density (solid curve) is the dominant factor in the unconditional density in Figure 1. For $\theta_2 = 0.1$, the situation is similar if slightly less pronounced. In case $\theta_2 = 0.75$, the solid and the dashed curve in Figure 2 get approximately equal weight, i.e., 0.51 and 0.49, respectively, resulting in a bimodal unconditional density in Figure 1. Finally, in case $\theta_2 = 1.2$, the weight of the solid curve is 0.12 while that of the dashed curve is 0.88; the resulting unconditional density in Figure 1 is unimodal but has a 'hump' in the left tail. For a detailed discussion of the conditional distributions and densities themselves, we refer to Section 3.3 of Leeb [1].

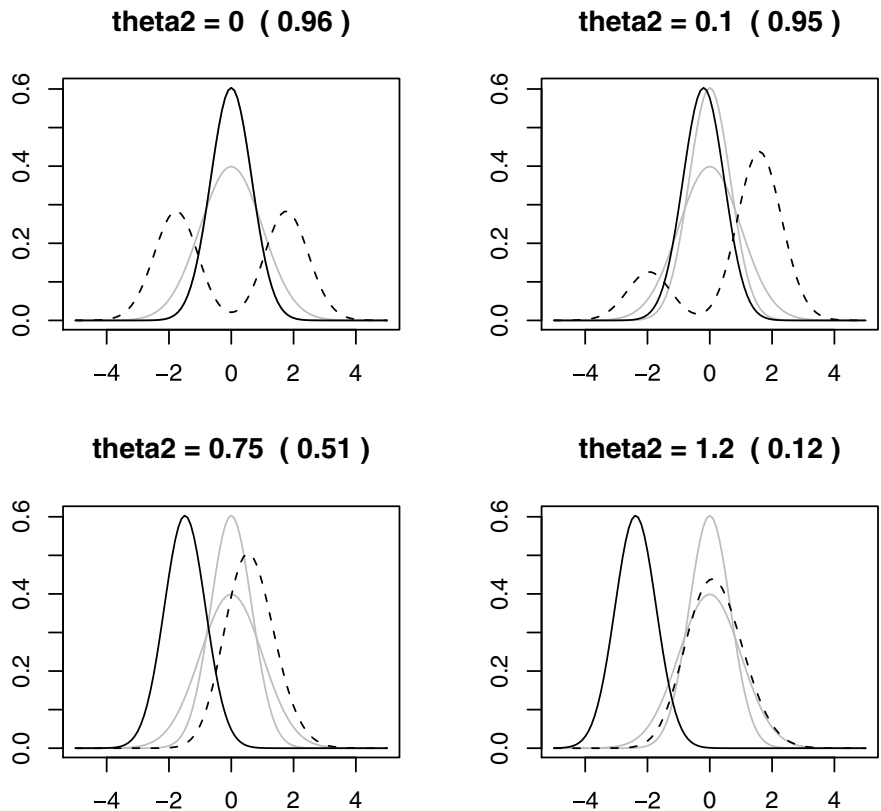

FIG 2. *The conditional density of $\sqrt{n}(\tilde{\theta}_1^* - \theta_1)$, conditional on selecting model $M_1$ (black solid line), and conditional on selecting model $M_2$ (black dashed line), for the same parameters as used for Figure 1. The number in parentheses in each panel header is the probability of selecting $M_1$, i.e., $\pi_{n,\theta,\sigma}^*(1)$. The gray curves are as in Figure 1.*

Results similar to Figure 1 and Figure 2 can be obtained for any other sample size (by appropriate choice of $\theta_2$ as noted above), and also for other choices of the critical value $c_2$ that is used by the model selectors. Larger values of $c_2$ result in model selectors that more strongly favor the smaller model $M_1$, and for which the phenomena observed above are more pronounced (see also Section 2.1 of Leeb and Pötscher [5] for results on the case where the critical value increases with sample size). Concerning the correlation coefficient $\rho$, we find that the shape of the conditional and of the unconditional densities is very strongly influenced by the magnitude of $|\rho|$, which we have chosen as $\rho = 0.75$ in figures 1 and 2 above. For larger values of $|\rho|$ we get similar but more pronounced phenomena. As $|\rho|$ gets smaller, however, these phenomena tend to be less pronounced. For example, if we plot the unconditional densities as in Figure 1 but with $\rho = 0.25$, we get four rather similar curves which altogether roughly resemble a Gaussian density except for some skewness. This is in line with the observation made in Section 3.3.1 that the unconditional distributions are Gaussian in the special case where $\tilde{\theta}_p(p)$ is uncorrelated with $A\tilde{\theta}(p)$ for each $p = \mathcal{O}+1, \ldots, P$. In the simple setting considered here, we have, in particular, that the distribution of $\sqrt{n}(\tilde{\theta}_1 - \theta_1)$ is Gaussian in the special case where $\rho = 0$.

## 4. An approximation result

In Theorem 4.2 below, we show that $G^*_{n,\theta,\sigma}(t)$ is close to $G_{n,\theta,\sigma}(t)$ in large samples, uniformly in the underlying parameters, where closeness is with respect to the total variation distance. (A similar result is provided in Leeb [1] for the conditional cdfs under slightly stronger assumptions.) Theorem 4.2 will be instrumental in the large-sample analysis in Section 5, because the large-sample behavior of $G^*_{n,\theta,\sigma}(t)$ is significantly easier to analyze. The total variation distance of two cdfs $G$ and $G^*$ on $\mathbb{R}^k$ will be denoted by $||G - G^*||_{TV}$ in the following. (Note that the relation $|G(t) - G^*(t)| \leq ||G - G^*||_{TV}$ always holds for each $t \in \mathbb{R}^k$. Thus, if $G$ and $G^*$ are close with respect to the total variation distance, then $G(t)$ is close to $G^*(t)$, uniformly in $t$. We shall use the total variation distance also for distribution functions $G$ and $G^*$ which are not necessarily normalized, i.e., in the case where $G$ and $G^*$ are the distribution functions of finite measures with total mass possibly different from one.)

Since the unconditional cdfs $G_{n,\theta,\sigma}(t)$ and $G^*_{n,\theta,\sigma}(t)$ can be linearly expanded in terms of $G_{n,\theta,\sigma}(t|p)\pi_{n,\theta,\sigma}(p)$ and $G^*_{n,\theta,\sigma}(t|p)\pi^*_{n,\theta,\sigma}(p)$, respectively, a key step for the results in this section is the following lemma.

**Lemma 4.1.** *For each $p$, $\mathcal{O} \leq p \leq P$, we have*

$$(4.1) \qquad \sup_{\substack{\theta \in \mathbb{R}^P \\ \sigma > 0}} \left|\left|G_{n,\theta,\sigma}(\cdot|p)\pi_{n,\theta,\sigma}(p) - G^*_{n,\theta,\sigma}(\cdot|p)\pi^*_{n,\theta,\sigma}(p)\right|\right|_{TV} \stackrel{n \to \infty}{\longrightarrow} 0.$$

This lemma immediately leads to the following result.

**Theorem 4.2.** *For the unconditional cdfs $G_{n,\theta,\sigma}(t)$ and $G^*_{n,\theta,\sigma}(t)$ we have*

(4.2) $$\sup_{\substack{\theta \in \mathbb{R}^P \\ \sigma > 0}} \|G_{n,\theta,\sigma} - G^*_{n,\theta,\sigma}\|_{TV} \overset{n \to \infty}{\longrightarrow} 0.$$

*Moreover, for each $p$ satisfying $\mathcal{O} \leq p \leq P$, the model selection probabilities $\pi_{n,\theta,\sigma}(p)$ and $\pi^*_{n,\theta,\sigma}(p)$ satisfy*

$$\sup_{\substack{\theta \in \mathbb{R}^P \\ \sigma > 0}} |\pi_{n,\theta,\sigma}(p) - \pi^*_{n,\theta,\sigma}(p)| \overset{n \to \infty}{\longrightarrow} 0.$$

By Theorem 4.2 we have, in particular, that

$$\sup_{\substack{\theta \in \mathbb{R}^P \\ \sigma > 0}} \sup_{t \in \mathbb{R}^k} |G_{n,\theta,\sigma}(t) - G^*_{n,\theta,\sigma}(t)| \overset{n \to \infty}{\longrightarrow} 0;$$

that is, the cdf $G_{n,\theta,\sigma}(t)$ is closely approximated by $G^*_{n,\theta,\sigma}(t)$ if $n$ is sufficiently large, uniformly in the argument $t$ and uniformly in the parameters $\theta$ and $\sigma$. The result in Theorem 4.2 does not depend on the scaling factor $\sqrt{n}$ and on the centering constant $A\theta$ that are used in the definitions of $G_{n,\theta,\sigma}(t)$ and $G^*_{n,\theta,\sigma}(t)$, cf. (2.4) and (2.5), respectively. In fact, that result continues to hold for arbitrary measurable transformations of $\tilde{\theta}$ and $\tilde{\theta}^*$. (See Corollary A.1 below for a precise formulation.)

Leeb [1] gives a result paralleling (4.2) for the conditional distributions of $A\tilde{\theta}$ and $A\tilde{\theta}^*$, conditional on the outcome of the model selection step. That result establishes closeness of the corresponding conditional cdfs uniformly not over the whole parameter space but over a slightly restricted set of parameters; cf. Theorem 4.1 in Leeb [1]. This restriction arose from the need to control the behavior of ratios of probabilities which vanish asymptotically. (Indeed, the probability of selecting the model of order $p$ converges to zero as $n \to \infty$ if the selected model is incorrect; cf. (5.6) below.) In the unconditional case considered in Theorem 4.2 above, this difficulty does not arise, allowing us to avoid this restriction.

## 5. Asymptotic results for the unconditional distributions and for the selection probabilities

We now analyze the large-sample limit behavior of $G_{n,\theta,\sigma}(t)$ and $G^*_{n,\theta,\sigma}(t)$, both in the fixed parameter case where $\theta$ and $\sigma$ are kept fixed while $n$ goes to infinity, and along sequences of parameters $\theta^{(n)}$ and $\sigma^{(n)}$. The main result in this section is Proposition 5.1 below. Inter alia, this result gives a complete characterization of all accumulation points of the unconditional cdfs (with respect to weak convergence) along sequences of parameters; cf. Remark 5.5. Our analysis also includes the model selection probabilities, as well as the case of local-alternative and fixed-parameter asymptotics.

The following conventions will be employed throughout this section: For $p$ satisfying $0 < p \leq P$, partition $Q = \lim_{n \to \infty} X'X/n$ as

$$Q = \begin{pmatrix} Q[p:p] & Q[p:\neg p] \\ Q[\neg p:p] & Q[\neg p:\neg p] \end{pmatrix},$$

where $Q[p:p]$ is a $p \times p$ matrix. Let $\Phi_{\infty,p}(t)$ be the cdf of a $k$-variate centered Gaussian random vector with covariance matrix $\sigma^2 A[p]Q[p:p]^{-1}A[p]'$, $0 < p \leq P$, and let $\Phi_{\infty,0}(t)$ denote the cdf of point-mass at zero in $\mathbb{R}^k$. Note that $\Phi_{\infty,p}(t)$ has a density with respect to Lebesgue measure on $\mathbb{R}^k$ if $p > 0$ and the matrix $A[p]$ has rank $k$; in this case, we denote the Lebesgue density of $\Phi_{\infty,p}(t)$ by $\phi_{\infty,p}(t)$. Finally, for $p = 1, \ldots, P$, define the quantities

$$\xi^2_{\infty,p} = (Q[p:p]^{-1})_{p,p},$$
$$\zeta^2_{\infty,p} = \xi^2_{\infty,p} - C_\infty^{(p)'}(A[p]Q[p:p]^{-1}A[p]')^- C_\infty^{(p)}, \text{ and}$$
$$b_{\infty,p} = C_\infty^{(p)'}(A[p]Q[p:p]^{-1}A[p]')^-,$$

where $C_\infty^{(p)} = A[p]Q[p:p]^{-1}e_p$, with $e_p$ denoting the $p$-th standard basis vector in $\mathbb{R}^p$. As the notation suggests, $\Phi_{\infty,p}(t)$ is the large-sample limit of $\Phi_{n,p}(t)$, $C_\infty^{(p)}$, $\xi_{\infty,p}$ and $\zeta_{\infty,p}$ are the limits of $C_n^{(p)}$, $\xi_{n,p}$ and $\zeta_{n,p}$, respectively, and $b_{n,p}z \to b_{\infty,p}z$ for each $z$ in the column-space of $A[p]$; cf. Lemma A.2 in Leeb [1]. With these conventions, we can characterize the large-sample limit behavior of the unconditional cdfs along sequences of parameters.

**Proposition 5.1.** *Consider sequences of parameters $\theta^{(n)} \in \mathbb{R}^P$ and $\sigma^{(n)} > 0$, such that $\sqrt{n}\theta^{(n)}$ converges to a limit $\psi \in (\mathbb{R} \cup \{-\infty, \infty\})^P$, and such that $\sigma^{(n)}$ converges to a (finite) limit $\sigma > 0$ as $n \to \infty$. Let $p_*$ denote the largest index $p$, $\mathcal{O} < p \leq P$, for which $|\psi_p| = \infty$, and set $p_* = \mathcal{O}$ if no such index exists. Then $G^*_{n,\theta^{(n)},\sigma^{(n)}}(t)$ and $G_{n,\theta^{(n)},\sigma^{(n)}}(t)$ both converge weakly to a limit cdf which is given by*

$$\Phi_{\infty,p_*}(t - A\delta^{(p_*)}) \prod_{q=p_*+1}^{P} \Delta_{\sigma\xi_{\infty,q}}(\delta^{(q)}_q + \psi_q, c_q\sigma\xi_{\infty,q})$$

(5.1) $$+ \sum_{p=p_*+1}^{P} \int_{z \leq t - A\delta^{(p)}} \left(1 - \Delta_{\sigma\zeta_{\infty,p}}(\delta^{(p)}_p + \psi_p + b_{\infty,p}z, c_p\sigma\xi_{\infty,p})\right) \Phi_{\infty,p}(dz)$$

$$\times \prod_{q=p+1}^{P} \Delta_{\sigma\xi_{\infty,q}}(\delta^{(q)}_q + \psi_q, c_q\sigma\xi_{\infty,q}),$$

*where*

(5.2) $$\delta^{(p)} = \begin{pmatrix} Q[p:p]^{-1}Q[p:\neg p] \\ -I_{P-p} \end{pmatrix} \psi[\neg p],$$

$p_* \leq p \leq P$ *(with the convention that $\delta^{(P)}$ is the zero-vector in $\mathbb{R}^P$ and, if necessary, that $\delta^{(0)} = -\psi$). Note that $\delta^{(p)}$ is the limit of the bias of $\tilde{\theta}(p)$ scaled by $\sqrt{n}$, i.e., $\delta^{(p)} = \lim_{n\to\infty} \sqrt{n}(\eta_n(p) - \theta^{(n)})$, with $\eta_n(p)$ given by (3.1) with $\theta^{(n)}$ replacing $\theta$; also note that $\delta^{(p)}$ is always finite, $p_* \leq p \leq P$.*

*The above statement continues to hold with convergence in total variation replacing weak convergence in the case where $p_* > 0$ and the matrix $A[p_*]$ has rank $k$, and in the case where $p_* < P$ and $\sqrt{n}A[\neg p_*]\theta^{(n)}[\neg p_*]$ is constant in $n$.*

**Remark 5.2.** Observe that the limit cdf in (5.1) is of a similar form as the finite-sample cdf $G^*_{n,\theta,\sigma}(t)$ as given in (3.7) (the only difference being that the right-hand side of (3.7) is the sum of $P - \mathcal{O} + 1$ terms while (5.1) is the sum of $P - p_* + 1$ terms, that quantities depending on the regressor matrix through $X'X/n$ in (3.7) are replaced by their corresponding limits in (5.1), and that the bias and mean of $\sqrt{n}\tilde{\theta}(p)$ in (3.7) are replaced by the appropriate large-sample limits in (5.1)). Therefore, the discussion of the finite-sample cdf $G^*_{n,\theta,\sigma}(t)$ given in Section 3.3 applies, mutatis mutandis, also to the limit cdf in (5.1). In particular, the cdf in (5.1) has a density with respect to Lebesgue measure on $\mathbb{R}^k$ if (and only if) $p_* > 0$ and $A[p_*]$ has rank $k$; in that case, this density can be obtained from (5.1) by differentiation. Moreover, we stress that the limit cdf is typically non-Gaussian. A notable exception where (5.1) reduces to the Gaussian cdf $\Phi_{\infty,P}(t)$ occurs in the special case where $\tilde{\theta}_q(q)$ and $A\tilde{\theta}(q)$ are asymptotically uncorrelated for each $q = p_* + 1, \ldots, P$.

Inspecting the proof of Proposition 5.1, we also obtain the large-sample limit behavior of the conditional cdfs weighted by the model selection probabilities, e.g., of $G_{n,\theta^{(n)},\sigma^{(n)}}(t|p)\pi_{n,\theta^{(n)},\sigma^{(n)}}(p)$ (weak convergence of not necessarily normalized cdfs $H_n$ to a not necessarily normalized cdf $H$ on $\mathbb{R}^k$ is defined as follows: $H_n(t)$ converges to $H(t)$ at each continuity point $t$ of the limit cdf, and $H_n(\mathbb{R}^k)$, i.e., the total mass of $H_n$ on $\mathbb{R}^k$, converges to $H(\mathbb{R}^k)$).

**Corollary 5.3.** *Assume that the assumptions of Proposition* 5.1 *are met, and fix $p$ with $\mathcal{O} \leq p \leq P$. In case $p = p_*$, $G_{n,\theta^{(n)},\sigma^{(n)}}(t|p_*)\pi_{n,\theta^{(n)},\sigma^{(n)}}(p_*)$ converges to the first term in* (5.1) *in the sense of weak convergence. If $p > p_*$, $G_{n,\theta^{(n)},\sigma^{(n)}}(t|p)\pi_{n,\theta^{(n)},\sigma^{(n)}}(p)$ converges weakly to the term with index $p$ in the sum in* (5.1). *Finally, if $p < p_*$, $G_{n,\theta^{(n)},\sigma^{(n)}}(t|p)\pi_{n,\theta^{(n)},\sigma^{(n)}}(p)$ converges to zero in total variation. The same applies to $G^*_{n,\theta^{(n)},\sigma^{(n)}}(t|p)\pi^*_{n,\theta^{(n)},\sigma^{(n)}}(p)$. Moreover, weak convergence can be strengthened to convergence in total variation in the case where $p > 0$ and $A[p]$ has rank $k$ (in that case, the weighted conditional cdf also has a Lebesgue density), and in the case where $p < P$ and $\sqrt{n}A[\neg p]\theta^{(n)}[\neg p]$ is constant in $n$.*

**Proposition 5.4.** *Under the assumptions of Proposition* 5.1, *the large-sample limit behavior of the model selection probabilities $\pi_{n,\theta^{(n)},\sigma^{(n)}}(p)$, $\mathcal{O} \leq p \leq P$, is as follows: For each $p$ satisfying $p_* < p \leq P$, $\pi_{n,\theta^{(n)},\sigma^{(n)}}(p)$ converges to*

$$(5.3) \quad (1 - \Delta_{\sigma\xi_{\infty,p}}(\delta^{(p)}_p + \psi_p, c_p\sigma\xi_{\infty,p})) \prod_{q=p+1}^{P} \Delta_{\sigma\xi_{\infty,q}}(\delta^{(q)}_q + \psi_q, c_q\sigma\xi_{\infty,q}).$$

*For $p = p_*$, $\pi_{n,\theta^{(n)},\sigma^{(n)}}(p_*)$ converges to*

$$(5.4) \quad \prod_{q=p_*+1}^{P} \Delta_{\sigma\xi_{\infty,q}}(\delta^{(q)}_q + \psi_q, c_q\sigma\xi_{\infty,q}).$$

*For each $p$ satisfying $\mathcal{O} \leq p < p_*$, $\pi_{n,\theta^{(n)},\sigma^{(n)}}(p)$ converges to zero. The above statements continue to hold with $\pi^*_{n,\theta^{(n)},\sigma^{(n)}}(p)$ replacing $\pi_{n,\theta^{(n)},\sigma^{(n)}}(p)$.*

**Remark 5.5.** With Propositions 5.1 and 5.4 we obtain a complete characterization of all possible accumulation points of the unconditional cdfs (with respect to weak convergence) and of the model selection probabilities, along arbitrary sequences of parameters $\theta^{(n)}$ and $\sigma^{(n)}$, provided that $\sigma^{(n)}$ is bounded away from zero and infinity: Let $\theta^{(n)}$ be any sequence in $\mathbb{R}^P$ and let $\sigma^{(n)}$ be a sequence satisfying $\sigma_* \leq \sigma^{(n)} \leq \sigma^*$ with $0 < \sigma_* \leq \sigma^* < \infty$. Since the set $(\mathbb{R} \cup \{-\infty, \infty\})^P$ as well as the set $[\sigma_*, \sigma^*]$ is compact, each subsequence contains a further subsequence for which the assumptions of Propositions 5.1 and 5.4 are satisfied. For example, each accumulation point of $G_{n,\theta^{(n)},\sigma^{(n)}}(t)$ (with respect to weak convergence) is of the form (5.1), where here $\psi$ and $\sigma$ are accumulation points of $\sqrt{n}\theta^{(n)}$ and $\sigma^{(n)}$, respectively (and where $p_*$ and the quantities $\delta^{(p)}$, $p_* \leq p \leq P$, are derived from $\psi$ as in Proposition 5.1). Of course, the same is true for $G^*_{n,\theta^{(n)},\sigma^{(n)}}(t)$. The same considerations apply, mutatis mutandis, to the weighted conditional cdfs considered in Corollary 5.3.

To study, say, the large-sample limit minimal coverage probability of confidence sets for $A\theta$ centered at $A\tilde\theta$, a description of all possible accumulation points of $G_{n,\theta^{(n)},\sigma^{(n)}}(t)$ with respect to weak convergence is useful; here $\theta^{(n)}$ can be any sequence in $\mathbb{R}^P$ and $\sigma^{(n)}$ can be any sequence bounded away from zero and infinity. In view of Remark 5.5, we see that each individual accumulation point can be reached along a particular sequence of regression parameters $\theta^{(n)}$, chosen such that the $\theta^{(n)}$ are within an $O(1/\sqrt{n})$ neighborhood of one of the models under consideration, say, $M_{p_*}$ for some $\mathcal{O} \le p_* \le P$. In particular, in order to describe all possible accumulation points of the unconditional cdf, it suffices to consider local alternatives to $\theta$.

**Corollary 5.6.** *Fix $\theta \in \mathbb{R}^P$ and consider local alternatives of the form $\theta + \gamma/\sqrt{n}$, where $\gamma \in \mathbb{R}^P$. Moreover, let $\sigma^{(n)}$ be a sequence of positive real numbers converging to a (finite) limit $\sigma > 0$. Then Propositions 5.1 and 5.4 apply with $\theta + \gamma/\sqrt{n}$ replacing $\theta^{(n)}$, where here $p_*$ equals $\max\{p_0(\theta), \mathcal{O}\}$ and $\psi[\neg p_*]$ equals $\gamma[\neg p_*]$ (in case $p_* < P$). In particular, $G^*_{n,\theta+\gamma/\sqrt{n},\sigma^{(n)}}(t)$ and $G_{n,\theta+\gamma/\sqrt{n},\sigma^{(n)}}(t)$ converge in total variation to the cdf in (5.1) with $p_* = \max\{p_0(\theta), \mathcal{O}\}$.*

In the case of fixed-parameter asymptotics, the large-sample limits of the model selection probabilities and of the unconditional cdfs take a particularly simple form. Fix $\theta \in \mathbb{R}^P$ and $\sigma > 0$. Clearly, $\sqrt{n}\theta$ converges to a limit $\psi$, whose $p_0(\theta)$-th component is infinite if $p_0(\theta) > 0$ (because the $p_0(\theta)$-th component of $\theta$ is non-zero in that case), and whose $p$-th component is zero for each $p > p_0(\theta)$. Therefore, Propositions 5.1 and 5.4 apply with $p_* = \max\{p_0(\theta), \mathcal{O}\}$, and either with $p_* < P$ and $\psi[\neg p_*] = (0, \ldots, 0)'$, or with $p_* = P$. In particular, $p_* = \max\{p_0(\theta), \mathcal{O}\}$ is the order of the smallest correct model for $\theta$ among the candidate models $M_\mathcal{O}, \ldots, M_P$. We hence obtain that $G^*_{n,\theta,\sigma}(t)$ and $G_{n,\theta,\sigma}(t)$ converge in total variation to the cdf

$$
(5.5) \quad \begin{aligned} &\Phi_{\infty,p_*}(t) \prod_{q=p_*+1}^{P} \Delta_{\sigma\xi_{\infty,q}}(0, c_q\sigma\xi_{\infty,q}) \\ &+ \sum_{p=p_*+1}^{P} \int_{z \le t} (1 - \Delta_{\sigma\zeta_{\infty,p}}(b_{\infty,p}z, c_p\sigma\xi_{\infty,p}))\Phi_{\infty,p}(dz) \\ &\times \prod_{q=p+1}^{P} \Delta_{\sigma\xi_{\infty,q}}(0, c_q\sigma\xi_{\infty,q}), \end{aligned}
$$

and the large-sample limit of the model selection probabilities $\pi_{n,\theta,\sigma}(p)$ and $\pi^*_{n,\theta,\sigma}(p)$ for $\mathcal{O} \le p \le P$ is given by

$$
(5.6) \quad \begin{cases} (1 - \Delta_{\sigma\xi_{\infty,p}}(0, c_p\sigma\xi_{\infty,p})) \prod_{q=p+1}^{P} \Delta_{\sigma\xi_{\infty,q}}(0, c_q\sigma\xi_{\infty,q}) & \text{if } p > p_*, \\ \prod_{q=p_*+1}^{P} \Delta_{\sigma\xi_{\infty,q}}(0, c_q\sigma\xi_{\infty,q}) & \text{if } p = p_*, \\ 0 & \text{if } p < p_* \end{cases}
$$

with $p_* = \max\{p_0(\theta), \mathcal{O}\}$.

**Remark 5.7.** (i) In defining the cdf $G_{n,\theta,\sigma}(t)$, the estimator has been centered at $\theta$ and scaled by $\sqrt{n}$; cf. (2.4). For the finite-sample results in Section 3, a different choice of centering constant (or scaling factor) of course only amounts to a translation (or rescaling) of the distribution and is hence inconsequential. Also, the results in Section 4 do not depend on the centering constant and on the scaling factor, because the total variation distance of two cdfs is invariant under a shift or rescaling of the argument. More generally, Lemma 4.1 and Theorem 4.2 extend to the distribution of arbitrary (measurable) functions of $\tilde{\theta}$ and $\tilde{\theta}^*$; cf. Corollary A.1 below.

(ii) We are next concerned with the question to which extent the limiting results given in the current section are affected by the choice of the centering constant. Let $d_{n,\theta,\sigma}$ denote a $P \times 1$ vector which may depend on $n$, $\theta$ and $\sigma$. Then centering at $d_{n,\theta,\sigma}$ leads to

$$(5.7) \quad P_{n,\theta,\sigma}\left(\sqrt{n}A(\tilde{\theta} - d_{n,\theta,\sigma}) \leq t\right) = G_{n,\theta,\sigma}\left(t + \sqrt{n}A(d_{n,\theta,\sigma} - \theta)\right).$$

The results obtained so far can now be used to describe the large-sample behavior of the cdf in (5.7). In particular, assuming that $\sqrt{n}A(d_{n,\theta,\sigma} - \theta)$ converges to a limit $\nu \in \mathbb{R}^k$, it is easy to verify that the large-sample limit of the cdf in (5.7) (in the sense of weak convergence) is given by the cdf in (5.5) with $t + \nu$ replacing $t$. If $\sqrt{n}A(d_{n,\theta,\sigma} - \theta)$ converges to a limit $\nu \in (\mathbb{R} \cup \{-\infty, \infty\})^k$ with some component of $\nu$ being either $\infty$ or $-\infty$, then the limit of (5.7) will be degenerate in the sense that at least one marginal distribution mass will have escaped to $\infty$ or $-\infty$. In other words, if $i$ is such that $|\nu_i| = \infty$, then the $i$-th component of $\sqrt{n}A(\tilde{\theta} - d_{n,\theta,\sigma})$ converges to $-\nu_i$ in probability as $n \to \infty$. The marginal of (5.7) corresponding to the finite components of $\nu$ converges weakly to the corresponding marginal of (5.5) with the appropriate components of $t + \nu$ replacing the appropriate components of $t$. This shows that, for an asymptotic analysis, any reasonable centering constant typically must be such that $Ad_{n,\theta,\sigma}$ coincides with $A\theta$ up to terms of order $O(1/\sqrt{n})$. If $\sqrt{n}A(d_{n,\theta,\sigma} - \theta)$ does not converge, accumulation points can be described by considering appropriate subsequences. The same considerations apply to the cdf $G^*_{n,\theta,\sigma}(t)$, and also to asymptotics along sequences of parameters $\theta^{(n)}$ and $\sigma^{(n)}$.

## Acknowledgments

I am thankful to Benedikt M. Pötscher for helpful remarks and discussions.

## Appendix A: Proofs for Section 4

*Proof of Lemma 4.1.* Consider first the case where $p > \mathcal{O}$. In that case, it is easy to see that $G_{n,\theta,\sigma}(t|p)\pi_{n,\theta,\sigma}(p)$ does not depend on the critical values $c_q$ for $q < p$ which are used by the model selection procedure $\hat{p}$ (cf. formula (3.9) above for $\pi_{n,\theta,\sigma}(p)$ and the expression for $G_{n,\theta,\sigma}(t|p)$ given in (16)–(18) of Leeb [1]). As a consequence, we conclude for $p > \mathcal{O}$ that $G_{n,\theta,\sigma}(t|p)\pi_{n,\theta,\sigma}(p)$ follows the same formula irrespective of whether $\mathcal{O} = 0$ or $\mathcal{O} > 0$. The same applies, mutatis mutandis, to $G^*_{n,\theta,\sigma}(t|p)\pi^*_{n,\theta,\sigma}(t)$. We hence may assume that $\mathcal{O} = 0$ in the following.

In the special case where $A$ is the $p \times P$ matrix $(I_p : 0)$ (which is to be interpreted as $I_P$ in case $p = P$), (4.1) follows from Lemma 5.1 of Leeb and Pötscher [3]. (In that result the conditional cdfs are such that the estimators are centered at $\eta_n(p)$ instead of $\theta$. However, this different centering constant does not affect the total variation distance; cf. Lemma A.5 in Leeb [1].) For the case of general $A$, write $\mu$ as shorthand for the conditional distribution of $\sqrt{n}(I_p : 0)(\tilde{\theta} - \theta)$ given $\hat{p} = p$ multiplied by $\pi_{n,\theta,\sigma}(p)$, $\mu^*$ as shorthand for the conditional distribution of $\sqrt{n}(I_p : 0)(\tilde{\theta}^* - \theta)$ given $\hat{p}^* = p$ multiplied by $\pi^*_{n,\theta,\sigma}(p)$, and let $\Psi$ denote the mapping $z \mapsto ((A[p]z)' : (-\sqrt{n}A[\neg p]\theta[\neg p])')'$ in case $p < P$ and $z \mapsto Az$ in case $p = P$. It is now easy to see that Lemma A.5 of Leeb [1] applies, and (4.1) follows.

It remains to show that (4.1) also holds with $\mathcal{O}$ replacing $p$. Having established (4.1) for $p > \mathcal{O}$, it also follows, for each $p = \mathcal{O}+1, \ldots, P$, that

$$\text{(A.1)} \qquad \sup_{\substack{\theta \in \mathbb{R}^P \\ \sigma > 0}} \left| \pi_{n,\theta,\sigma}(p) - \pi^*_{n,\theta,\sigma}(p) \right| \stackrel{n \to \infty}{\longrightarrow} 0,$$

because the modulus in (A.1) is bounded by

$$||G_{n,\theta,\sigma}(\cdot|p)\pi_{n,\theta,\sigma}(p) - G^*_{n,\theta,\sigma}(\cdot|p)\pi^*_{n,\theta,\sigma}(p)||_{TV}.$$

Since the model selection probabilities sum up to one, we have $\pi_{n,\theta,\sigma}(\mathcal{O}) = 1 - \sum_{p=\mathcal{O}+1}^{P} \pi_{n,\theta,\sigma}(p)$, and a similar expansion holds for $\pi^*_{n,\theta,\sigma}(\mathcal{O})$. By this and the triangle inequality, we see that (A.1) also holds with $\mathcal{O}$ replacing $p$. Now (4.1) with $\mathcal{O}$ replacing $p$ follows immediately, because the conditional cdfs $G_{n,\theta,\sigma}(t|\mathcal{O})$ and $G^*_{n,\theta,\sigma}(t|\mathcal{O})$ are both equal to $\Phi_{n,\mathcal{O}}(t - \sqrt{n}A(\eta_n(\mathcal{O}) - \theta))$, cf. (10) and (14) of Leeb [1], which is of course bounded by one. $\square$

*Proof of Theorem 4.2.* Relation (4.2) follows from Lemma 4.1 by expanding $G^*_{n,\theta,\sigma}(t)$ as in (3.4), by expanding $G_{n,\theta,\sigma}(t)$ in a similar fashion, and by applying the triangle inequality. The statement concerning the model selection probabilities has already been established in the course of the proof of Lemma 4.1; cf. (A.1) and the attending discussion. $\square$

**Corollary A.1.** *For each $n$, $\theta$ and $\sigma$, let $\Psi_{n,\theta,\sigma}(\cdot)$ be a measurable function on $\mathbb{R}^P$. Moreover, let $R_{n,\theta,\sigma}(\cdot)$ denote the distribution of $\Psi_{n,\theta,\sigma}(\tilde{\theta})$, and let $R^*_{n,\theta,\sigma}(\cdot)$ denote the distribution of $\Psi_{n,\theta,\sigma}(\tilde{\theta}^*)$. (That is, say, $R_{n,\theta,\sigma}(\cdot)$ is the probability measure induced by $\Psi_{n,\theta,\sigma}(\tilde{\theta})$ under $P_{n,\theta,\sigma}(\cdot)$.) We then have*

$$\text{(A.2)} \qquad \sup_{\substack{\theta \in \mathbb{R}^P \\ \sigma > 0}} \left| \left| R_{n,\theta,\sigma}(\cdot) - R^*_{n,\theta,\sigma}(\cdot) \right| \right|_{TV} \stackrel{n \to \infty}{\longrightarrow} 0.$$

*Moreover, if $R_{n,\theta,\sigma}(\cdot|p)$ and $R^*_{n,\theta,\sigma}(\cdot|p)$ denote the distributions of $\Psi_{n,\theta,\sigma}(\tilde{\theta})$ conditional on $\hat{p} = p$ and of $\Psi_{n,\theta,\sigma}(\tilde{\theta}^*)$ conditional on $\hat{p}^* = p$, respectively, then*

$$\text{(A.3)} \qquad \sup_{\substack{\theta \in \mathbb{R}^P \\ \sigma > 0}} \left| \left| R_{n,\theta,\sigma}(\cdot|p)\pi_{n,\theta,\sigma}(p) - R^*_{n,\theta,\sigma}(\cdot|p)\pi^*_{n,\theta,\sigma}(p) \right| \right|_{TV} \stackrel{n \to \infty}{\longrightarrow} 0.$$

*Proof.* Observe that the total variation distance of two cdfs is unaffected by a change of scale or a shift of the argument. Using Theorem 4.2 with $A = I_P$, we hence obtain that (A.2) holds if $\Psi_{n,\theta,\sigma}$ is the identity map. From this, the general case follows immediately in view of Lemma A.5 of Leeb [1]. In a similar fashion, (A.3) follows from Lemma 4.1. $\square$

## Appendix B: Proofs for Section 5

Under the assumptions of Proposition 5.1, we make the following preliminary observation: For $p \geq p_*$, consider the scaled bias of $\tilde{\theta}(p)$, i.e., $\sqrt{n}(\eta_n(p) - \theta^{(n)})$, where $\eta_n(p)$ is defined as in (3.1) with $\theta^{(n)}$ replacing $\theta$. It is easy to see that

$$\sqrt{n}(\eta_n(p) - \theta^{(n)}) = \begin{pmatrix} (X[p]'X[p])^{-1}X[p]'X[\neg p] \\ -I_{P-p} \end{pmatrix} \sqrt{n}\theta^{(n)}[\neg p],$$

where the expression on the right-hand side is to be interpreted as $\sqrt{n}\theta^{(n)}$ and as the zero vector in $\mathbb{R}^P$ in the cases $p = 0$ and $p = P$, respectively. For $p$ satisfying $p_* \leq p < P$, note that $\sqrt{n}\theta^{(n)}[\neg p]$ converges to $\psi[\neg p]$ by assumption, and that this limit is finite by choice of $p \geq p_*$. It hence follows that $\sqrt{n}(\eta_n(p) - \theta^{(n)})$ converges to the limit $\delta^{(p)}$ given in (5.2). From this, we also see that $\sqrt{n}\eta_{n,p}(p)$ converges to $\delta_p^{(p)} + \psi_p$, which is finite for each $p > p_*$; for $p = p_*$, this limit is infinite in case $|\psi_{p_*}| = \infty$. Note that the case where the limit of $\sqrt{n}\eta_{n,p_*}(p_*)$ is finite can only occur if $p_* = \mathcal{O}$. It will now be convenient to prove Proposition 5.4 first.

*Proof of Proposition 5.4.* In view of Theorem 4.2, it suffices to consider $\pi^*_{n,\theta^{(n)},\sigma^{(n)}}(p)$. This model selection probability can be expanded as in (3.5)–(3.6) with $\theta^{(n)}$ and $\sigma^{(n)}$ replacing $\theta$ and $\sigma$, respectively. Consider first the individual $\Delta$-functions occurring in these formulas, i.e.,

(B.1) $$\Delta_{\sigma^{(n)}\xi_{n,q}}(\sqrt{n}\eta_{n,q}(q), c_q\sigma^{(n)}\xi_{n,q}),$$

$\mathcal{O} < q \leq P$. For $q > p_*$, recall that $\sqrt{n}\eta_{n,q}(q)$ converges to the finite limit $\delta_q^{(q)} + \psi_q$ as we have seen above, and it is elementary to verify that the expression in (B.1) converges to $\Delta_{\sigma\xi_{\infty,q}}(\delta_q^{(q)} + \psi_q, c_q\sigma\xi_{\infty,q})$. For $q = p_*$ and $p_* > \mathcal{O}$, we have seen that the limit of $\sqrt{n}\eta_{n,p_*}(p_*)$ is infinite, and it is easy to see that (B.1) with $p_*$ replacing $q$ converges to zero in this case.

From the above considerations, it immediately follows that $\pi^*_{n,\theta^{(n)},\sigma^{(n)}}(p)$ converges to the limit in (5.3) if $p > p_*$, and to the limit in (5.4) if $p = p_*$. To show that $\pi^*_{n,\theta^{(n)},\sigma^{(n)}}(p)$ converges to zero in case $p$ satisfies $\mathcal{O} \leq p < p_*$, it suffices to observe that here $\pi^*_{n,\theta^{(n)},\sigma^{(n)}}(p)$ is bounded by the expression in (B.1) with $p_*$ replacing $q$. As we have seen above, $\sqrt{n}|\eta_{n,p_*}(p_*)|$ converges to infinity, such that this upper bound converges to zero as $n \to \infty$. □

*Proof of Proposition 5.1.* Again, it suffices to consider $G^*_{n,\theta^{(n)},\sigma^{(n)}}(t)$ in view of Theorem 4.2. Recall that this cdf can be written as in (3.4). We first consider the individual terms $G^*_{n,\theta^{(n)},\sigma^{(n)}}(t|p)\pi^*_{n,\theta^{(n)},\sigma^{(n)}}(p)$ for $p = \mathcal{O}, \ldots, P$. In case $p$ satisfies $\mathcal{O} \leq p < p_*$, note that $\pi^*_{n,\theta^{(n)},\sigma^{(n)}}(p) \to 0$ by Proposition 5.4. Hence, $G^*_{n,\theta^{(n)},\sigma^{(n)}}(t|p)\pi^*_{n,\theta^{(n)},\sigma^{(n)}}(p)$ converges to zero in total variation.

In the remaining cases, i.e., for $p$ satisfying $p_* \leq p \leq P$, it is elementary to verify that Proposition 5.1 of Leeb [1] applies to $G^*_{n,\theta^{(n)},\sigma^{(n)}}(t|p)$, where the quantity $\beta$ in that paper equals $A\delta^{(p)}$ in our setting. In particular, that result gives the limit of the conditional cdf in the sense of weak convergence (because $\delta^{(p)}$ is finite). Consider first the case $p > p_*$. From Proposition 5.4, we obtain the limit of $\pi^*_{n,\theta^{(n)},\sigma^{(n)}}(p)$. Combining the resulting limit expression with the limit expression for $G^*_{n,\theta^{(n)},\sigma^{(n)}}(t|p)$ as obtained by Proposition 5.1 of Leeb [1], we see that

$G^*_{n,\theta^{(n)},\sigma^{(n)}}(t|p)\pi^*_{n,\theta^{(n)},\sigma^{(n)}}(p)$ converges weakly to

$$\text{(B.2)} \quad \int_{\substack{z \in \mathbb{R}^k \\ z \leq t - A\delta^{(p)}}} \left(1 - \Delta_{\sigma\zeta_{\infty,p}}(\delta^{(p)}_p + \psi_p + b_{\infty,p}z, c_p\sigma\xi_{\infty,p})\right) \Phi_{\infty,p}(dz)$$
$$\times \prod_{q=p+1}^{P} \Delta_{\sigma\xi_{\infty,q}}(\delta^{(q)}_q + \psi_q, c_q\sigma\xi_{\infty,q}).$$

In case $p = p_*$ and $p_* > \mathcal{O}$, we again use Proposition 5.1 of Leeb [1] and Proposition 5.4 to obtain that the weak limit of $G^*_{n,\theta^{(n)},\sigma^{(n)}}(t|p_*)\pi^*_{n,\theta^{(n)},\sigma^{(n)}}(p_*)$ is of the form (B.2) with $p_*$ replacing $p$. Since $|\psi_{p_*}|$ is infinite, the integrand in (B.2) reduces to one, i.e., the limit is given by

$$\Phi_{\infty,p_*}(t - A\delta^{(p_*)}) \prod_{q=p_*+1}^{P} \Delta_{\sigma\xi_{\infty,q}}(\delta^{(q)}_q + \psi_q, c_q\sigma\xi_{\infty,q}).$$

Finally, consider the case $p = p_*$ and $p_* = \mathcal{O}$. Arguing as above, we see that $G^*_{n,\theta^{(n)},\sigma^{(n)}}(t|\mathcal{O})\pi^*_{n,\theta^{(n)},\sigma^{(n)}}(\mathcal{O})$ converges weakly to

$$\Phi_{\infty,\mathcal{O}}(t - A\delta^{(\mathcal{O})}) \prod_{q=\mathcal{O}+1}^{P} \Delta_{\sigma\xi_{\infty,q}}(\delta^{(q)}_q + \psi_q, c_q\sigma\xi_{\infty,q}).$$

Because the individual model selection probabilities $\pi^*_{n,\theta^{(n)},\sigma^{(n)}}(p)$, $\mathcal{O} \leq p \leq P$, sum up to one, the same is true for their large-sample limits. In particular, note that (5.1) is a convex combination of cdfs, and that all the weights in the convex combination are positive. From this, we obtain that $G^*_{n,\theta^{(n)},\sigma^{(n)}}(t)$ converges to the expression in (5.1) at each continuity point $t$ of the limit expression, i.e., $G^*_{n,\theta^{(n)},\sigma^{(n)}}(t)$ converges weakly. (Note that a convex combination of cdfs on $\mathbb{R}^k$ is continuous at a point $t$ if each individual cdf is continuous at $t$; the converse is also true, provided that all the weights in the convex combination are positive.) To establish that weak convergence can be strengthened to convergence in total variation under the conditions given in Proposition 5.1, it suffices to note, under these conditions, that $G^*_{n,\theta^{(n)},\sigma^{(n)}}(t|p)$, $p_* \leq p \leq P$, converges not only weakly but also in total variation in view of Proposition 5.1 of Leeb [1]. □